\documentclass[
final
]{dmtcs-episciences}

\usepackage[utf8]{inputenc}
\usepackage{subfigure}

\usepackage[round]{natbib}
\usepackage{amssymb,amsmath}

\usepackage{appendix}
\newcommand{\bea}{\begin{eqnarray}}
\newcommand{\ena}{\end{eqnarray}}
\newcommand{\beas}{\begin{eqnarray*}}
\newcommand{\enas}{\end{eqnarray*}}
\newcommand{\beq}{\begin{equation}}
\newcommand{\enq}{\end{equation}}
\def\qed{\hfill \mbox{\rule{0.5em}{0.5em}}}

\newcommand{\ignore}[1]{}

\newtheorem{theorem}{Theorem}[section]
\newtheorem{corollary}{Corollary}[section]

\newtheorem{proposition}{Proposition}[section]

\newtheorem{remark}{Remark}[section]

\author[Recep Altar Çiçeksiz et. al]{Recep Altar Çiçeksiz \affiliationmark{1}  \and Yunus Emre Demirci \affiliationmark{2}  \and \"{U}m\.{i}t I\c{s}lak \affiliationmark{3}}
\title{The variance  and the asymptotic distribution of the length of longest $k$-alternating subsequences}
\affiliation{
  Ume\aa{} University, Department of Mathematics, Ume\aa, Sweden\\
  Queen's University, Department of Mathematics and Statistics, Kingston, Ontario, Canada\\
  Bo\u{g}azi\c{c}i University, Department of Mathematics, Istanbul, Turkey}
\begin{document}
% This is only used if you are compiling for a volume before vol 25
% \publicationdetails{VOL}{2015}{ISS}{NUM}{SUBM}
% This is the new form of collecting the data, starting with vol 25
\publicationdata{vol. 25:1}{2023}{12}{10.46298/dmtcs.10296}{2022-11-12; None}{2023-03-12}

\maketitle
\begin{abstract}
We obtain an explicit formula  for the variance of the number of $k$-peaks  in a uniformly random permutation. This is then used to obtain an asymptotic formula for the variance of  the  length of  longest $k$-alternating subsequence in     random permutations. Also a central limit is proved for the latter statistic.
\end{abstract}

\section{Introduction}
\label{sec:in}

Letting
$(a_i)_{i=1}^n$ be a sequence of real numbers, a
subsequence $a_{i_k}$, where $1 \leq i_1 < \ldots< i_k \leq n $, is
called an \emph{alternating subsequence} if
$a_{i_1}>a_{i_2}<a_{i_3}> \cdots.$ The \emph{length of the
longest alternating subsequence of} $(a_i)_{i=1}^n$ is defined to be the largest integer $q$ such that  $(a_i)_{i=1}^n$  has an alternating subsequence of
length $q$. %We denote this by $\operatorname{as}_{n,1}$.  
Denoting  the symmetric group on $n$ letters by $S_n$, an   alternating subsequence of a permutation $\sigma \in S_n$ refers to an alternating subsequence corresponding to the sequence $\sigma(1),\sigma(2), \ldots, \sigma(n)$.
 See \cite{stanley2008} for a survey on the topic.

The purpose of this manuscript is to study a generalization of the  length of longest alternating subsequences in uniformly random permutations. 
  Letting  $\sigma \in S_n$, a subsequence $1 \leq i_1 < i_2 < \ldots < i_t \leq n$ is said to be    $k$-\emph{alternating} for $\sigma$ if 
\[\sigma(i_1) \geq \sigma(i_2)+k, \quad \sigma(i_2)+k \leq \sigma(i_3), \quad \sigma(i_3) \geq \sigma(i_4)+k, \cdots.\]
In other words,  the subsequence is $k$-alternating if it is alternating and additionally
\[|\sigma(i_j)-\sigma(i_{j+1})| \geq k, \quad j \in [t-1],\] where we set $[m] = \{1,\ldots,m\}$ for $m \in \mathbb{N}$. Below the length of the longest $k$-alternating subsequence of $\sigma \in S_n$ is denoted by $\operatorname{as}_{n,k}(\sigma)$, or simply $\operatorname{as}_{n,k}$.

%For example, letting $\mathbf{a}=(3,1,7,4,2,6,5),$ one can easily see  that   $\operatorname{as}_{7,1}(\mathbf{a})=6$ and $\operatorname{as}_{7,2}(\mathbf{a})=5$. 
  
%  When we talk about alternating subsequences of a permutation $\sigma  \in S_n$, the definitions are meant to be for the sequence $a_j = \sigma(j)$. A general reference  on the longest alternating subsequence problem is  \cite{stanley2008}. 

Let us also define $k$-peaks and $k$-valleys which will be intermediary tools   to understand  the longest $k$-alternating subsequences. Let $\sigma=\sigma(1) \ldots \sigma(n) \in S_n$. We say that a section  $\sigma(i) \ldots \sigma(j)$ of the permutation $\sigma$ is a $k$\emph{-up} ($k$\emph{-down}, resp.) if $i <j$ and $\sigma(j) -\sigma(i) \geq k$ ($\sigma(i)-\sigma(j) \geq k$, resp.). We say that the section is $k$\emph{-ascending} if it satisfies: 
\begin{itemize}
    \item $\sigma(i)=\min\{\sigma(i), \ldots, \sigma(j)\}$ and $\sigma(j)=\max\{\sigma(i), \ldots,\sigma(j)\}$, and
    \item the section $\sigma(i) \ldots \sigma(j)$ is a $k$-up, and 
    \item there is no $k$-down in $\sigma(i) \ldots \sigma(j)$, i.e. for any $i \leq s<t \leq j$, we have $\sigma(s)-\sigma(t)<k$.
\end{itemize}
If also there is no $k$-ascending section that contains $\sigma(i) \ldots \sigma (j)$, it is called a \emph{maximal $k$-ascending section}. In this case,  $\sigma(i)$,   $\sigma(j)$ are called    a $k$\emph{-valley} and  a  $k$\emph{-peak} of $\sigma$, respectively. 
 
A maximal  $k$-descending section $\sigma(i)\ldots \sigma(j)$ can be defined similarly, and this time  $\sigma(i)$,   $\sigma(j)$ are called a  $k$-peak and a $k$-valley of $\sigma$, respectively. An alternative description can be given  as  in \cite{cai}.
\begin{proposition}\label{propn:char}
Let $\sigma=\sigma(1) \sigma(2) \ldots \sigma(n) \in S_n$, $i \in [n]$ and $1 \leq k \leq n-1$. Then $\sigma(i)$ is a $k$-peak if and only if it satisfies both of the following two properties: \par
(i) If there is an $s>i$ with $\sigma(s)>\sigma(i)$, then there is a $k$-down $\sigma(i)\ldots \sigma(j)$ in $\sigma(i) \ldots \sigma(s)$. \par
(ii) If there is an $s<i$ with $\sigma(s)>\sigma(i)$, then there is a $k$-up $\sigma(j) \ldots \sigma(i)$ in $\sigma(s) \ldots \sigma(i)$. 
\end{proposition}

Considering the case where $\sigma $ is a uniformly random permutation, our purpose in present paper is to study $\operatorname{Var}(\operatorname{as}_{n,k})$ and to show that $\operatorname{as}_{n,k}$ satisfies a central limit theorem. The statistic $\operatorname{Var}(\operatorname{as}_{n,k})$ is well understood for the case $k=1$.  Indeed, Stanley  proved in \cite{stanley2008} that 
$$ \mathbb{E}[\operatorname{as}_{n,1}] = \frac{4 n+1}{6} \quad  \text{and}  \quad 
\operatorname{Var}[ \operatorname{as}_{n,1} ]=\frac{8 n}{45}-\frac{13}{180}.$$
It was later shown in  \cite{houdre}  and  \cite{romik}     that $\operatorname{as}_{n,1}$ satisfies a central limit theorem, and convergence rates for the normal approximation were obtained in    \cite{I:2018}.  All these limiting distribution results rely on the simple fact that $ \operatorname{as}_{n,1}$ can be represented as a sum of $m$-dependent random variables (namely, the indicators of local extrema) and they then  use the well-established theory of such sequences.

Regarding the general $k$, Armstrong conjectured in \cite{armstrong}  that $\mathbb{E} [\operatorname{as}_{n,k}]=\frac{4(n-k)+5}{6}$.   Pak and
  Pemantle \cite{pak2014longest} then used probabilistic methods to prove that $\mathbb{E} [\operatorname{as}_{n,k}]$ is asymptotically $\frac{2(n-k)}{3}+O\left(n^{2 / 3}\right)$.  
  
Let us very briefly mention their approach. For $x \in(0,1)$,
a vector $\mathbf{y}=\left(y_{1}, \ldots, y_{n}\right) \in[0,1]^{n}$ is said to be  $x$\emph{-alternating} if $(-1)^{j+1}\left(y_{j}-y_{j+1}\right) \geqslant x$ for all $1 \leqslant j \leqslant n-1$. Given a vector $\mathbf{y}=\left(y_{1}, \ldots, y_{n}\right) \in[0,1]^{n}$, a subsequence $1 \leq i_1 < i_2 < \ldots < i_r \leq n$ is said to be $x$\emph{-alternating for} $\mathbf{y}$ if
\[| y_{i_j}-y_{i_{j+1}}| \geq x, \quad j \in [r-1].\]
Denoting the length of the longest $x-$alternating subsequence of a random vector $\mathbf{y}$, with Lebesgue measure on $[0,1]^n$ as its distribution, by $\operatorname{as_{n,x}}(\mathbf{y})$, their main observation was: 
If  $Z$ is a binomial random variable with parameters $n$ and $1-x$, then  
$$
\operatorname{as}_{n, x}(\mathbf{y}) \stackrel{\mathcal{D}}{=} \operatorname{as}_{Z,1}
$$
(Here, $\stackrel{\mathcal{D}}{=}$ means equality in distribution). That is, they concluded that $\operatorname{as}_{n, x}(\mathbf{y})$ has the same distribution as the length of the longest  ordinary alternating subsequence of a random permutation on $S_Z$.   Afterwards, using 
$ \mathbb{E}[\operatorname{as}_{n,1}] = \frac{4 n+1}{6}$ and $\operatorname{Var}(\operatorname{as}_{n,1} )=\frac{8 n}{45}-\frac{13}{180}$, they proved 
$$\mathbb{E} [\operatorname{as}_{n, x}] =\frac{2}{3} n(1-x)+\frac{1}{6} \quad  \text{and} \quad \operatorname{Var}( \operatorname{as}_{n, x} )=(1-x)(2+5 x) \frac{4 n}{45}.$$  Further, for suitable $x_1$ and $x_2$, they showed that $\mathbb{E} [\operatorname{as}_{n, x_{2}}] \leqslant \mathbb{E} [\operatorname{as}_{n, k}] \leqslant \mathbb{E}[\operatorname{as}_{n, x_{1}}]$ and  in this way they are  able to bound $\mathbb{E} [\operatorname{as}_{n, k}]$.

A closely related problem to the longest alternating subsequence problem is that of calculating the longest zigzagging subsequence. For a given permutation $\sigma$, denoting its vertical flip by $\tilde{\sigma}$,   a subsequence is said to be \emph{zigzagging} if it is alternating for either $\sigma$ or $\tilde{\sigma}$. The   $k$-zigzagging case is defined similarly.   We will be using the notation  $\operatorname{zs}_{n,k}$ for the length of the longest $k$-zigzagging subsequence in the sequel. Note  that in exactly half of the permutations, $\operatorname{as}_{n,k}$ and $\operatorname{zs}_{n,k}$ are   equal to each other, and in the other half the length of the $k$-zigzagging subsequence is exactly one more than the length of the $k$-alternating subsequence. This is seen via the involution map $I: \sigma(1)\sigma(2) \ldots \sigma(n) \to (n+1-\sigma(1))(n+1-\sigma(2))\ldots(n+1-\sigma(n))$ as noted in \cite{cai}.
Therefore 
\begin{equation}\label{atozconn}
\mathbb{E}[\operatorname{z s}_{k}]=\mathbb{E}[\operatorname{a s}_{k}]+1/2.
\end{equation}
Cai proved  in 2015  that $\mathbb{E}[\operatorname{z s}_{k}] =\frac{2(n-k)+4}{3}$, and then combining this with \eqref{atozconn}, solved the   Armstrong conjecture  \cite{cai}.

Our first  result in this paper is  an asymptotic  formula for $\operatorname{Var}(\operatorname{as}_{n, k})$. Namely, we will prove $$\operatorname{Var}(\operatorname{as}_{n, k})=\dfrac{8(n-k)}{45}+ O (\sqrt{n}).$$ In order to obtain this   result,  we first study the number of $k$-peaks $P$ in random permutations and show that $$\operatorname{Var} (P) = \dfrac{2(n-k) + 4}{45}.$$ 

Our second result is a central limit theorem for $\operatorname{as}_{n, k}$:
\[\dfrac{\operatorname{as}_{n,k}-\mathbb{E}[\operatorname{as}_{n,k}]}{\sqrt{\operatorname{Var}(\operatorname{as}_{n,k})}} \longrightarrow_{d} \mathcal{G},\]
where $\mathcal{G}$ is the standard normal distribution and where $\rightarrow_d$ is used for convergence in distribution.

The rest of the paper is organized as follows. Next section proves our formulas for the variances of $P$ and $\operatorname{as}_{n, k}$.  In Section \ref{sec:clt}, we prove the central limit theorem for $\operatorname{as}_{n, k}$. 

% You may scarsely use \clearpage to advance to a new page if this
% improves the readability of the document structure
\clearpage

\section{The variances of $P$ and $\operatorname{as}_{n, k}$}

Next result gives an exact formula for the variance of the number of   $k$-peaks $P$ in a uniformly random permutation.  
%First we are interested in the number of   $k$-peaks $P$ in a uniformly random permutation.  

\begin{theorem}\label{thm:varpeaks}
Let $P$ be the number of   $k$-peaks  in a uniformly random permutation in $S_n$.  We have  $$\operatorname{Var} (P) = \dfrac{2(n-k) + 4}{45}. $$
\end{theorem}

We will prove 
Theorem \ref{thm:varpeaks} after providing a corollary related to the length of longest $k$-alternating subsequence of  a uniformly random permutation.  Note that we have $\operatorname{as}_{n, k} =  2 P + E $ where $|E| \leq 1$ for any $n,k$. Thus, $\operatorname{Var} (\operatorname{as}_{n, k} )= 4 \operatorname{Var} (P)+ \operatorname{Var} (E) + 2 \operatorname{Cov} (P,E).$ Here, clearly $\operatorname{Var} (E) \leq 1$ and by Cauchy-Schwarz inequality $|\operatorname{Cov} (P,E)| \leq 2 \sqrt{\operatorname{Var}(P)} \sqrt{\operatorname{Var}(E)} \leq C_0 \sqrt{n}$ where $C_0$ is a constant independent of $n$ and $k$.
We now obtain the following.
 \begin{corollary}\label{thm:main}
Let  $\operatorname{as}_{n, k}$ be the length of longest $k$-alternating subsequence of  a uniformly random permutation  in $S_n$. Then,   $$\operatorname{Var}(\operatorname{as}_{n, k})=\dfrac{8(n-k)}{45}+ O (\sqrt{n}).$$ In particular, when $k = o(n)$, $\operatorname{Var}(\operatorname{as}_{n, k}) \sim \frac{8n}{45}$ as $n \rightarrow \infty$.
\end{corollary}
\begin{remark} In setting of Corollary \ref{thm:main}, we conjecture that 
$\operatorname{Var}(\operatorname{as}_{n, k})=\dfrac{8(n-k)}{45}+\dfrac{19}{180}$.  Although we have a heuristic derivation of this equality,  we were not able to justify it rigorously. 
\end{remark}
 Now, let us proceed to the proof of Theorem \ref{thm:varpeaks}.

\vspace{0.1in}

\textit{Proof of Theorem \ref{thm:varpeaks}.}
%We are willing to compute $\operatorname{Var} (\operatorname{as}_{n,k}) = \mathbb{E}[\operatorname{as}_{n,k}^2] - (\mathbb{E}[\operatorname{as}_{n,k}])^2$. 
 Below    $P_i$ is  the indicator of $i$ being a $k$-peak\footnote{Note that when we say $i$ is a $k$-peak,  we consider $i$ to be an element in the image of the permutation, not an element of the domain of the permutation. If the position $i$ is considered in domain of the permutation, we will be emphasizing it there.}, i.e. 
$$P_i:=\begin{cases}
1, & i \;\; \text{is a k-peak}, \\
0, & \text{otherwise}.
\end{cases}$$
In particular,  $$P=\sum_{i=1}^n P_i.$$ 
We are willing to compute  
$$ \operatorname{Var} (P) =  \operatorname{Var}\left(\sum_{i=1}^n P_i \right) =   \mathbb{E} \left[ \left(\sum_{i=1}^n P_i \right)^2  \right] - \left( \mathbb{E} \left[  \sum_{i=1}^n P_i  \right]\right)^2. $$
Recall from \cite{cai} that  
\begin{equation}\label{EP}
 \mathbb{E} \left[  \sum_{i=1}^n P_i  \right] = \mathbb{E}[P] =  \frac{1}{2}  \mathbb{E}[\operatorname{z s}_{k}] =\frac{n-k+2}{3}.
\end{equation}

Let us next analyze   
$$\mathbb{E}\left[\left( \sum_{i=1}^n P_i\right)^2 \right] = \sum_{i=1}^n \mathbb{E}[P_i^2]+2\sum_{i<j}\mathbb{E}[P_iP_j].$$
Denoting  the probability that $i$ is a $k$-peak  by $p_{n,k}(i)$ and  the probability that both $i,j$ are $k$-peaks by $p_{n,k}(i,j)$, we may rewrite this last equation as   $$\mathbb{E}\left[\left( \sum_{i=1}^n P_i\right)^2 \right] = \sum_{i=1}^n p_{n,k}(i)+2\sum_{i<j}p_{n,k}(i,j).$$
We already know from \eqref{EP} that the first sum on the right-hand side  is $\frac{n-k+2}{3}$.   We  are then left with   calculating  $p_{n,k}(i,j)$.

With the definition of $k$-peaks in mind, for given $i$ and $j$, we can divide $[n]\setminus \{i\}$ and $[n]\setminus \{j\}$ into three sets according to the following partitions respectively. The first partition is with respect to $i$:
\begin{align*}
A_i=&\lbrace \ell: 1 \leq  \ell \leq i-k \rbrace, \\
B_i=&\lbrace  \ell: i-k+1 \leq  \ell \leq i-1 \rbrace, \\
C_i=&\lbrace  \ell: i+1 \leq  \ell \leq n \rbrace,
\end{align*}
and the second partition is with respect to $j$:
\begin{align*}
A_j=&\lbrace  \ell: 1 \leq  \ell \leq j-k \rbrace, \\
B_j=&\lbrace  \ell: j-k+1 \leq  \ell \leq j-1 \rbrace, \\
C_j=&\lbrace  \ell: j+1 \leq  \ell \leq n \rbrace.
\end{align*}
Assuming without loss of generality that $i<j$, observe 
\begin{align*}
i<j& \implies A_i \subset A_j \\
i<j& \implies C_j \subset C_i.
\end{align*}

By Proposition \ref{propn:char}, we observe that for $i$ to be a $k$-peak, there should be at least one element from $A_i$ between any element of $C_i$ and $i$, and similarly for $j$ to be a $k$-peak, there should be at least one element from $A_j$ between any element of $C_j$ and $j$. To ensure these two properties, we will place the elements accordingly.

Our procedure for placing the elements starts with placing $A_i \cup \{i\}$ in a row $a_1a_2 \ldots a_{i-k+1}$ arbitrarily. 
Leaving the insertion of  the elements in $A_j \setminus A_i$ to the end of the argument, we will next focus on placing the elements of  $C_i$ and $C_j$. %We begin with the following observation. 
%Since, $A_i \subset A_j$, the elements in $A_i$ affect $j$ to be a $k$-peak as well as $i$. By this we mean the fact that if we only look at $j$ being a $k$-peak, we would have put all the elements in $A_j$ and $j$ in row, yet since the elements in $A_j \setminus A_i$ do not affect $i$ being a $k$-peak, and if these elements are inserted before an element of $C_i$, then $i$ would not be a $k$-peak, we put them after we deal with the elements of $C_i$ and $C_j$.   
%Before moving on to placing the elements that require more care, an observation is required.  
%By the inclusions   stated above, we see that the elements in $C_i \cap C_j=C_j$ affect both $i$ and $j$ being a $k$-peak,   and the elements in $C_i \setminus C_j$ affect only $j$ being a $k$-peak. 
Note that by the observation in previous paragraph, in order to have  $i$ and $j$ as $k$-peaks,   the two places next to $i$ are not available for the elements in $C_i \setminus C_j$, and the four places next to $i$ and $j$ are not available for the elements in $C_i \cap C_j=C_j$.

Now, let us focus on the elements of $C_i \setminus C_j=\{i+1,\ldots,j\}$. There are $|A_i \cup \{i\}|=i-k+1$ elements that are placed in a row. Thus, we have $i-k+2$ vacant spots
 for the element $i-k+2$ to be inserted into the row $a_1a_2\ldots a_{i-k+1}$. Since the two places next to $i$ are prohibited, we see that 
$$\mathbb{P}(\{i+1\} \; \text{does not prevent } i,j \text{ being a } k\text{-peak})=\frac{i-k}{i-k+2}.$$
Now, we have $i+k+3$
 vacant spots for the element $i+2$, and the two places next to $i$ are prohibited, and so,
$$\mathbb{P}(\{i+2\} \; \text{does not prevent }  i,j \text{ being a } k\text{-peak})=\frac{i-k+1}{i-k+3}.$$
Continuing in this manner, we see that when we arrive at $j$, which is the last element to be inserted in from the set $C_i \setminus C_j$, we have $i-k+(j-i+1)=j-k+1$ many vacant places, and the two places next to $i$ are prohibited, and then
$$\mathbb{P}(\{j\} \; \text{does not prevent } i,j \text{ being a } k\text{-peak})=\frac{j-k-1}{j-k+1}.$$
More generally, for $t=1,\ldots,j- i$, we have $$\mathbb{P}\left(\{i+t\} \; \text{does not prevent  } i,j \text{  being a } k\text{-peak}\right) = \dfrac{i-k+t-1}{i-k+t+1}.$$
Therefore, 
\begin{small}
\begin{eqnarray*}
 \mathbb{P}\left(C_i \setminus C_j \; \text{ does not prevent } i,j    \text{ being a } k\text{-peak}\right) 
&=&\mathbb{P}\left(\bigcap_{t=1}^{j-i}\{i+t\} \; \text{does not prevent  } i,j \text{  being a } k\text{-peak}\right)\\
%&=&\prod_{t=1}^{j-i}\mathbb{P}\left(\{i+t\} \; \text{does not prevent } i,j \text{ being a } k\text{-peak}\right)\\
&=&\prod_{t=1}^{j-i}\dfrac{i-k+t-1}{i-k+t+1} \\
&=& \dfrac{(i-k)(i-k+1)}{(j-k)(j-k+1)}.    
\end{eqnarray*}
\end{small}

Now, let us focus on the elements of $C_i \cap C_j=C_j=\{j+1,\ldots,n\}$. Recall that there are four prohibited places  for these elements to be inserted. We have $j-k+2$ many vacant places to insert $j+1$ into but four of these are prohibited. Thus, 
$$\mathbb{P}(\{j+1\} \; \text{does not prevent } i,j \text{ being a } k\text{-peak})=\frac{j-k-2}{j-k+2}.$$
Similar to the analysis in $C_i \setminus C_j$, continuing in this manner, we have $n=j+(n-j)$, and  in the end we will have $j-k+(n-j+1)=n-k+1$ many vacant places to insert $n$, and four of these are prohibited. So,
$$\mathbb{P}(n \; \text{does not prevent }  i,j \text{ being a } k\text{-peak})=\frac{n-k-3}{n-k+1}.$$
We may generalize this to obtain  
$$\mathbb{P}\left( \{j+t\} \; \text{does not prevent } i,j \text{ being a } k\text{-peak}\right) = \dfrac{j-k+t-3}{j-k+t+1} $$ for $t = 1,\ldots,n-j$.  We then obtain 
\begin{small}
\begin{eqnarray*}
 \mathbb{P}( C_i \cap C_j \; \text{does not prevent }  i,j \text{ being a } k\text{-peak} )
&=&\mathbb{P}\left( \bigcap_{t=1}^{n-j}\{j+t\} \; \text{does not prevent } i,j \text{ being a } k\text{-peak}\right)\\
%&=&\prod_{t=1}^{n-j}\mathbb{P}(\{j+t\} \; \text{does not prevent } i,j \text{ being a } k\text{-peak})\\
&=&\prod_{t=1}^{n-j}\dfrac{j-k+t-3}{j-k+t+1} \\
&=&\dfrac{(j-k-2)(j-k-1)(j-k)(j-k+1)}{(n-k-2)(n-k-1)(n-k)(n-k+1)}.    
\end{eqnarray*}
\end{small}

Note that we can multiply the probabilities (here, and above in the case of $C_i \setminus C_j$),  since in essence what we are doing is conditioning on the event that the previous added elements do not prevent $i,j$ being a $k$-peak. Now, clearly, the elements of $A_j \setminus A_i$ are in $B_i \cup C_i$. Since the elements that are in $C_i$ have been inserted, we will then be done  once we  insert  the elements of $B_i$ and $B_j$. But the elements in the sets $B_i$ and $B_j$ have no effect on $i$ and $j$ being a $k$-peak (once the elements from $C_i$ and $C_j$ are placed), and so  we may insert them in any place. Thus, overall, we have 
\begin{eqnarray*}
p_{n,k}(i,j)&=& \mathbb{P}( \text{the set}\; C_i \cup C_j \; \text{does not prevent } i,j \text{ being a } k\text{-peak} )\\
&=& \mathbb{P} ( C_i \setminus C_j \; \text{does not prevent } i,j \text{ being a } k\text{-peak} ) \\
&& \times \mathbb{P} ( C_i \cap C_j \; \text{does not prevent } i,j \text{ being a } k\text{-peak} )\\
&=& \dfrac{(i-k)(i-k+1)}{(j-k)(j-k+1)}\dfrac{(j-k-2)(j-k-1)(j-k)(j-k+1)}{(n-k-2)(n-k-1)(n-k)(n-k+1)}\\
&=& \dfrac{(i-k)(i-k+1)(j-k-2)(j-k-1)}{(n-k-2)(n-k-1)(n-k)(n-k+1)}.    
\end{eqnarray*}
These add up to 
\begin{eqnarray*}
\sum_{i<j}p_{n,k}(i,j) &=& \sum_{i=k+1}^n\sum_{j=i+1}^n\dfrac{(i-k)(i-k+1)(j-k-2)(j-k-1)}{(n-k-2)(n-k-1)(n-k)(n-k+1)} \\
&=& \dfrac{1}{90}(5k-5n+3)(k-n-2),
\end{eqnarray*}
where the sum is computed fairly easily noting that essentially we are summing the consecutive integers and squares of consecutive integers. 
Therefore  we obtain
\begin{eqnarray*}
\mathbb{E}\left[\left(\sum_{i=1}^n P_i \right)^2 \right] &=&   \sum_{i=1}^n p_{n,k}(i) +2\sum_{i<j} p_{n,k}(i,j)  
 = \dfrac{n-k+2}{3}+\dfrac{1}{45}(5k-5n+3)(k-n-2) \\ &=& \dfrac{n-k+2}{3}\left(1+\dfrac{1}{15}(5n-5k-3)\right)
=\dfrac{1}{45}(n-k+2)(5n-5k+12).
\end{eqnarray*}

Using this we arrive at 
\begin{eqnarray*}
\operatorname{Var}(P)   &=& \mathbb{E} \left[ \left(\sum_{i=1}^n P_i \right)^2  \right] - \left( \mathbb{E} \left[  \sum_{i=1}^n P_i  \right]\right)^2  \\  
&=& \dfrac{1}{45}(n-k+2)(5n-5k+12)-  \left(\dfrac{n-k+2}{3} \right)^2  = \dfrac{2(n-k) + 4}{45} 
\end{eqnarray*}
as asserted in Theorem \ref{thm:varpeaks}. \hfill $\square$

\section{A Central Limit theorem}\label{sec:clt}
In this section, we will prove the following central limit theorem.
\begin{theorem}
Let $k$ be a fixed  positive integer.   Then the length of the longest $k$-alternating subsequence $\operatorname{as}_{n,k}$ of a uniformly random permutation satisfies a central limit theorem,
\[\dfrac{\operatorname{as}_{n,k}-\mathbb{E}[\operatorname{as}_{n,k}]}{\sqrt{\operatorname{Var}(\operatorname{as}_{n,k})}} \longrightarrow_{d} \mathcal{G},\]
where $\mathcal{G}$ is the standard normal distribution.
\end{theorem}

The proof involves a suitable truncation argument that allows us to reduce the problem to proving a central limit theorem for sums of locally dependent random variables for which a theory is already available.  Since the length of the longest $k$ alternating sequence differs from  twice the number of $k$ peaks by at most 1, we may focus on the number of peaks.   For any $i$, let $P_i$ be the random variable that is $1$ if the value $i$ is a $k$-peak and zero otherwise as before. Also recall $P = P_1 + \cdots + P_n$. We know that $P_i = 1$ precisely when 
\begin{itemize}
    \item Scanning to the right of the value $i$, we encounter an element in $[i-k]$ before we encounter an element in $[i+1, n]$. It is permitted that we do not encounter an element from $[i+1, n]$ at all. 
    \item Scanning to the left of the value $i$, we encounter an element in $[i-k]$ before we encounter an element in $[i+1, n]$. It is permitted that we do not encounter an element from $[i+1, n]$ at all. 
\end{itemize}

Our approach to getting a central limit theorem is to define a suitable truncation that can be computed using local data. There are a number of theorems that establish central limit behaviour for variables with only local correlations and this approach has been employed in a number of situations. 

Note that the condition on $P_i = 1$ can be restated as 
\begin{itemize}
    \item There is an index $j > \sigma^{-1}(i)$ such that $i-k \geq \sigma(j)$ and such that 
    \[i= \operatorname{max}_{s \in [\sigma^{-1}(i), j]} \sigma(s), \quad \sigma(j)= \operatorname{min}_{s \in [\sigma^{-1}(i), j]} \sigma(s).\]
    \item There is an index $j < \sigma^{-1}(i)$ such that $i-k \geq \sigma(j)$ and such that 
    \[i= \operatorname{max}_{s \in [j,\sigma^{-1}(i)]} \sigma(s), \quad \sigma(j)= \operatorname{min}_{s \in [j, \sigma^{-1}(i)]} \sigma(s).\]
\end{itemize}
Note that we might need to scan far to the left and right in order to determine whether a value is a $k$-peak or not and thus we will have long range dependence. We will show that ignoring long range interactions does not change the statistic very much. 

Fix a number $m$ that we will specify later. Let $Y_i = 1$ if we can determine that $i$ is a $k$-peak by only looking at $m$ positions to the left and right of $i$. Precisely, let $Y_i = 1$ if 
\begin{itemize}
    \item There is an index $j \in  [\sigma^{-1}(i), \sigma^{-1}(i)+m]$ such that  $i-k \geq \sigma(j)$ and such that 
    \[i= \operatorname{max}_{s \in [\sigma^{-1}(i), j]} \sigma(s), \quad \sigma(j)= \operatorname{min}_{s \in [\sigma^{-1}(i), j]} \sigma(s).\]
    \item There is an index $j \in  [\sigma^{-1}(i)-m, \sigma^{-1}(i)]$ such that  $i-k \geq \sigma(j)$ and such that 
    \[i= \operatorname{max}_{s \in [j,\sigma^{-1}(i)]} \sigma(s), \quad \sigma(j)= \operatorname{min}_{s \in [j, \sigma^{-1}(i)]} \sigma(s).\]
\end{itemize}

If $Y_i = 1$, we call it a \emph{local $k$-peak} (suppressing the reference to $m$). Note that any local $k$-peak is a $k$-peak and thus, $Y_i \leq P_i$. We should next understand the case where  $Y_i = 0$ and $P_i = 1$. Note that if $i \leq k$, then $P_i = Y_i = 0$. 

If $\sigma^{-1}(i) \in [m+1]$, there is no issue when scanning to the left. However, if we scan to the right and this event happens, then the $m$ indices to the right should have values in $[i-k+1, i-1]$. The probability of this is at most $\left(\frac{k-1}{n-1}\right)^m$. Similarly, the probability of this event when $\sigma^{-1}(i) \in [n-m, n]$ is at most $\left(\frac{k-1}{n-1}\right)^m$.

If $\sigma^{-1}(i) \in [m+2, n-m-2]$, the event can only happen if the $2m$ positions, $m$ to the left and $m$ to the right take values in $[i-k-1, i-1]$ and the probability of this is at most $\left(\frac{k-1}{n-1} \right)^{2m}$. 

Putting these together, recalling $Y_i \leq P_i$, and  denoting  the total variation distance by $d_{TV}$, we see that 
\begin{align*}
    d_{TV}(P_i, Y_i) &=\frac{1}{2}\sum_{j=0}^1 |\mathbb{P}(P_i=j)-\mathbb{P}(Y_i=j)| \\
    &=\frac{1}{2}(\mathbb{P}(Y_i=0)-\mathbb{P}(P_i=0)+\mathbb{P}[P_i=1]-\mathbb{P}(Y_i=1))   \\
    &=\frac{1}{2}(2(\mathbb{P}(P_i=1)-\mathbb{P}(Y_i=1))) \\
    &= \mathbb{P}(Y_i = 0, \, P_i = 1) \\ 
    &\leq \dfrac{2m+2}{n} \left(\dfrac{k-1}{n-1}\right)^m + \dfrac{n-2m-2}{n}\left(\dfrac{k-1}{n-1}\right)^{2m}.
\end{align*}

This implies  
\begin{align}\label{eqn:tvbound} 
\nonumber d_{TV}(P_1 + \ldots + P_n,  Y_1 + \ldots + Y_n) 
%\leq& \sum_{i=1}^n \mathbb{P}(P_i \neq Y_i) \\
 \leq& \dfrac{(2m+2)(n)}{n} \left(\dfrac{k-1}{n-1}\right)^m \\ 
\nonumber & + \dfrac{(n-2m-2)(n)}{n}\left(\dfrac{k-1}{n-1}\right)^{2m},\\
\nonumber \leq& (2m+2)\left(\dfrac{k}{n}\right)^m+n\,\left(\dfrac{k}{n}\right)^{2m},\\
\leq& 3n \left(\dfrac{k}{n}\right)^m.
\end{align}

When   $k$ is  fixed,  taking $m = 3$ suffices for our purpose.  Note in particular that
\begin{align}\label{bound}
\mathbb{P}(Y_1 + \ldots + Y_n < P_1 + \cdots + P_n) = o\left(\dfrac{1}{n}\right),
\end{align}
when $m$ is chosen appropriately. 

Next we  will show that $Y = Y_1 + \cdots + Y_n$ satisfies a central limit theorem. Let $Z_i = 1$ if the position $i$ is a local $k$-peak and $0$ otherwise. It is immediate that $Z = Z_1 + \cdots + Z_n$ and $Y$ have the same distribution. We let $Z$ be such a random variable for which $(P,Z)$ and $(P,Y)$ have the same distribution. Further, note that  the variables $Z_i$ have the property that $Z_i$ and $Z_j$ are independent if $|i - j| > 2m$. 

There are a number of related theorems that guarantee central limit behaviour for sums of locally dependent variables. A result due to Rinott \cite{rinott} will suffice for our purpose. The version we give is a slight variation of the one discussed in \cite{raic}. 

\begin{theorem} 
Let $U_1, \ldots, U_n$ be random variables such that $U_i$ and $U_j$ are independent when $|i-j| > 2m$. Setting $U = U_1 + \cdots + U_n$, we have  \[d_{K}\left(\dfrac{U-\mathbb{E}[U]}{\sqrt{\operatorname{Var}(U)}}, \mathcal{G} \right) \leq C(2m+1)\sqrt{\dfrac{\sum_{i = 1}^{n} \mathbb{E} |U_i|^3}{(\operatorname{Var}(U))^{3/2}}},\]  where $d_K$ is the Kolmogorov distance. 
\end{theorem}

We will now apply this result for $Z = Z_1 + \cdots + Z_n$. For this purpose
we need a lower bound on the variance of the random variable $Z$. Recall that the variance of $P$ is $\Omega(n)$ and let us show that the same holds for $Z$. 

We have   
\begin{align*} 
\sqrt{\operatorname{Var}(Z)} \geq& \sqrt{\operatorname{Var}(P)} - \sqrt{\operatorname{Var}(P-Z)}\\ 
\geq& \sqrt{\operatorname{Var}(P)} - \sqrt{\mathbb{E}[(P- Z)^2]} \\
\geq& \sqrt{\operatorname{Var}(P)} - \sqrt{n} \sqrt{\mathbb{E} |P - Z|}\\
\geq& \sqrt{\operatorname{Var}(P)} - \sqrt{n} \sqrt{\mathbb{E} [  |P - Z | \mid P \neq Z ] \mathbb{P}(P \neq Z)}\\
\geq& \sqrt{\operatorname{Var}(P)} - \sqrt{n} \sqrt{n}  \sqrt{o \left(\frac{1}{n} \right)} 
\tag*{using (\ref{bound})}\\
=& \Omega(\sqrt{n}) - o(\sqrt{n})\\
=&  \Omega(\sqrt{n}).
\end{align*}

Also, the $Z_i$ are Bernoulli random variables and thus $\sum_{i = 1}^{n} \mathbb{E} |Z_i|^3 = O(n)$. This shows that 
\[d_{K}\left(\dfrac{Z-\mathbb{E}[Z]}{\sqrt{\operatorname{Var}(Z)}},\mathcal{G} \right) \leq O\left(\dfrac{m}{n^{1/4}}\right),\]
proving  that  when $k$ is fixed, we have a central limit theorem,
\[\dfrac{Z-\mathbb{E}[Z]}{\sqrt{\operatorname{Var}(Z)}} \longrightarrow_{d} \mathcal{G}.\]
Together with the total variation distance bound between $P$ and $Z$, and noting that convergence in $d_{TV}$ implies convergence in $d_K$, we conclude that $P$ satisfies a central limit theorem. Since   $\operatorname{as}_{n,k}$ differs from $2P$ by at most 1, the same holds for it as well after proper centering and scaling.

\begin{remark}
The arguments given in this section carry over to certain cases where $k$ grows with $n$. For example, considering the case $k= \gamma n$ for constant $\gamma$,  the quantity $3n\left(\frac{k}{n}\right)^m$ in \eqref{eqn:tvbound} can be made $o(1/n)$ by choosing $m$ suitably. To see this, letting $\alpha>1$,   suppose $\frac{1}{n^\alpha}=3n\left(\frac{k}{n}\right)^m$. Since $\gamma=\frac{k}{n}$,  we then have $$n^{-1-\alpha}=3(\gamma)^m,$$ and then 
$m=\frac{(-1-\alpha)\operatorname{log}(n)-\operatorname{log}(3))}{\operatorname{log}(\gamma)}$. We can choose $\alpha=2$ so that $m =  \dfrac{-3\operatorname{log}(n)}{\operatorname{log}(\gamma)}$. Note that $m>0$ since $\log\left(\frac{k}{n}\right)<0.$
\end{remark}

\begin{remark}
In notation of the Introduction, if we were to prove a central limit theorem for $\operatorname{as_{n,x}}$, then that would be straightforward. This is thanks to the fact that it can be written as a random sum (where the number of summands is binomial) of locally dependent variables, and that central limit theorem for such cases are already available. See, for example, \cite{randomsums}.
\end{remark}
\bigskip

\acknowledgements We would like to thank Mohan Ravichandran for helpful discussions, especially towards the the local $k$-peaks argument used in the proof of the central limit theorem.  Third author is supported partially by BAP grant 20B06P

\nocite{*}
%\bibliographystyle{abbrvnat}
% use the following instead if you encounter problems 
\bibliographystyle{alpha}
\bibliography{sample-dmtcs}

\begin{thebibliography}{10}
\providecommand{\natexlab}[1]{#1}
\providecommand{\url}[1]{\texttt{#1}}
\expandafter\ifx\csname urlstyle\endcsname\relax
  \providecommand{\doi}[1]{doi: #1}\else
  \providecommand{\doi}{doi: \begingroup \urlstyle{rm}\Url}\fi

\bibitem[Armstrong(2014)]{armstrong}
D.~Armstrong.
\newblock Enumerative combinatorics problem session.
\newblock \emph{in Oberwolfach Report No}, 12 2014.

\bibitem[Cai(2015)]{cai}
T.~W. Cai.
\newblock Average length of the longest k-alternating subsequence.
\newblock \emph{Journal of Combinatorial Theory, Series A}, 134:\penalty0
  51--57, 2015.

\bibitem[Houdr\'e and Restrepo(2010)]{houdre}
C.~Houdr\'e and R.~Restrepo.
\newblock A probabilistic approach to the asymptotics of the length of the
  longest alternating subsequence.
\newblock \emph{Electronic Journal of Combinatorics}, 17, 2010.

\bibitem[Islak(2016)]{randomsums}
U.~Islak.
\newblock Asymptotic results for random sums of dependent random variables.
\newblock \emph{Statistics and Probability Letters}, 109:\penalty0 22--29,
  2016.

\bibitem[Islak(2018)]{I:2018}
U.~Islak.
\newblock Descent-inversion statistics in riffle shuffles.
\newblock \emph{Turkish Journal of Mathematics}, 42\penalty0 (2):\penalty0
  502--514, 2018.

\bibitem[Pak and Pemantle(2015)]{pak2014longest}
I.~Pak and R.~Pemantle.
\newblock On the longest k-alternating subsequence.
\newblock \emph{Electronic Journal of Combinatorics}, 22\penalty0 (1), 2015.

\bibitem[Raic(2003)]{raic}
M.~Raic.
\newblock Normal approximation by stein's method.
\newblock \emph{In Proceedings of the 7th Young Statisticians Meeting}, pages
  71--97, 2003.

\bibitem[Rinott(1994)]{rinott}
Y.~Rinott.
\newblock On normal approximation rates for certain sums of dependent random
  variables.
\newblock \emph{Journal of Computational and Applied Mathematics}, 55\penalty0
  (2):\penalty0 135--143, 1994.

\bibitem[Romik(2011)]{romik}
D.~Romik.
\newblock Local extrema in random permutations and the structure of longest
  alternating subsequences.
\newblock \emph{Discrete Mathematics and Theoretical Computer Science}, 2011.

\bibitem[Stanley(2008)]{stanley2008}
R.~Stanley.
\newblock Longest alternating subsequences of permutations.
\newblock \emph{Michigan Mathematical Journal}, 57:\penalty0 675--687, 2008.

\end{thebibliography}
\label{sec:biblio}

\end{document}